\documentclass[a4paper,12pt]{article}
\usepackage{comment}
\usepackage{cite}
\usepackage{amsmath}
\usepackage{amssymb}
\usepackage{amsfonts}
\usepackage[T1]{fontenc}
\usepackage[utf8]{inputenc}
\usepackage{graphicx}
\usepackage{fancyhdr}
\usepackage{float}
\usepackage{xcolor}
\usepackage{authblk}
\usepackage{mathrsfs}
\usepackage{empheq}
\usepackage[hyphens]{url}
\usepackage{hyperref} 
\usepackage[]{breakurl}

\usepackage{geometry}
\begin{document}

\title{An integral representation for $\zeta(4)$}

\author[$\dagger$]{Jean-Christophe {\sc Pain}$^{1,2,}$\footnote{jean-christophe.pain@cea.fr}\\
\small
$^1$CEA, DAM, DIF, F-91297 Arpajon, France\\
$^2$Universit\'e Paris-Saclay, CEA, Laboratoire Mati\`ere en Conditions Extr\^emes,\\ 
91680 Bruy\`eres-le-Ch\^atel, France
}

\maketitle

\begin{abstract}
In this note, we propose an integral representation for $\zeta(4)$, where $\zeta$ is the Riemann zeta function. The corresponding expression is obtained using relations for polylogarithms. A possible generalization to any even argument of the zeta function is considered.
\end{abstract}

\section{Introduction}

The Riemann zeta function for even values of the argument is given by
\begin{equation}
\zeta(2k)=\sum_{n=1}^{\infty}\frac1{n^{2k}} \ = \ \frac{|B_{2k}| \ (2\pi)^{2k}} {2 \, (2k)!},     
\end{equation}
where $B_n$ are Bernoulli numbers. A number of integral representations of the Riemann zeta function do exit. For instance, in the case where $\Re(s) > 1$, the following representation is known since Euler:
\begin{equation}
\zeta(s)={\frac {1}{\Gamma (s)}}\int _{0}^{1}{\frac {(-\ln u)^{s-1}}{1-u}}\,du,
\end{equation}
where $\Gamma$ represents the usual Gamma function. By making a change of variables, the latter expression turns out to be equivalent to
\begin{equation}
\zeta(s)=\frac{1}{\Gamma(s)}\int_0^{\infty}\frac{t^{s-1}}{e^t-1}\,dt.
\end{equation}
The latter formula can be interpreted as a particular case of a general transformation of Dirichlet series. It means that $s \longmapsto \zeta(s)\Gamma(s)$ is the Mellin transform of the function $t\mapsto 1/(e^t-1)$. Alternatively, using the Euler-Mac-Laurin formula\cite{Tenenbaum2008,Cartier1992}, one obtains
\begin{eqnarray}
\zeta (s)&=&{1 \over {s-1}}+{1 \over 2}+\sum _{k=2}^{n}B_{k}{\frac {s(s+1)\dots (s+k-2)}{k!}}\nonumber\\
& &-{\frac {s(s+1)\dots (s+n-1)}{n!}}\int _{1}^{\infty }B_{n}(x-[x])x^{-s-n}\,dx,
\end{eqnarray}
where $B_n(x)$ are the usual Bernoulli polynomials. Using contour integration, one finds also
\begin{equation}
\zeta(s)=\frac{\mathrm e^{-\mathrm i\pi s}\Gamma(1-s)}{2\mathrm i\pi}\oint_C{\frac{u^{s-1}}{\mathrm e^u-1}\,du},
\end{equation}
where $\mathscr{C}$ represents a contour following the real axis, embedding 0 and run from $-\infty$ to $+\infty$ counterclockwise. Borwein and Borwein obtained the following intriguig relation \cite{Borwein1995}:
\begin{equation}
\zeta(4)=\frac{2}{11\pi}\int_0^{\pi}\theta^2\ln^2\left[2\cos\left(\frac{\theta}{2}\right)\right]d\theta.
\end{equation}
The well-known value of $\zeta(4)$ can be obtained in different ways. For instance, let us introduce the $2\pi$-periodic function $f$ defined, $\forall x\in \left[-\pi ,\pi \right]$, by $f(x)=x^{2}$ and let us calculate its complex Fourier coefficients
\begin{equation}
c_{0}=\frac{1}{2\pi}\int_{-\pi}^{\pi}t^{2}\,dt={\frac {\pi ^{2}}{3}} 
\end{equation}
and
\begin{equation}
c_{n}=\frac{1}{2\pi}\int_{-\pi }^{\pi }t^{2}e^{-int}\,dt=2\frac{(-1)^n}{n^2}.
\end{equation}
The Fourier series of $f$ thus reads
\begin{equation}
f(x)=c_0+\sum_{n\in\mathbb{Z}^*}c_n\,e^{inx}.
\end{equation}
Using the Parseval-Plancherel formula:
\begin{equation}
\sum_{n=-\infty}^{\infty}|c_n|^2=\frac{1}{2\pi}\int_0^{2\pi}f^2(x)\,dx,
\end{equation}
we get
\begin{equation}
\frac {1}{2\pi}\int_{-\pi }^{\pi }t^{4}\,dt={\frac 
{\pi^{4}}{9}}+8\sum_{n=1}^{\infty}{\frac{1}{n^{4}}} 
\end{equation}
and thus
\begin{equation}
\sum_{n=1}^{\infty}{\frac {1}{n^{4}}}=\frac{\pi^{4}}{90}=\frac{1}{8}\left(\frac{\pi ^{4}}{5}-\frac{\pi^{4}}{9}\right).
\end{equation}
In the next section, we propose an integral relation for $\zeta(4)$. To our knowledge, the formula was not published elsewhere. It consists of an integral similar to the ``Fermi-Dirac'' type ones encountered in statistical physics.

\section{New integral relation}

Considering the integral
\begin{equation}
\mathscr{I}=\int_{-\infty}^{+\infty} z^2\, \frac{\ln(1+e^z)}{1+e^z}\,dz,
\end{equation}
let us start with the change of variables $u=e^z$, yielding
\begin{equation}
\mathscr{I}=\int_0^{\infty}\frac{(\ln u)^2\ln(1+u)}{u(1+u)}\,du.
\end{equation}
Let us now make the change of variable: $u=t/(1-t)$. One has $t=u/(1+u)$ and
\begin{equation}\label{lo}
\mathscr{I}=-\int_0^{1}\frac{\left[\ln (t)-\ln(1-t)\right]^2\ln(1-t)}{t}\,dt,  
\end{equation}
where
\begin{equation}
\mathrm{Li}_{s}(z)=\sum _{k=1}^{\infty }{z^{k} \over k^{s}}
\end{equation}
is the usual polylogarithm ($\mathrm{Li}_2$ is usually referred to as the ``dilogarithm''). One has in particular
\begin{equation}
\mathrm{Li}_{n}(1)=\zeta(n),
\end{equation}
where $\zeta$ is the Riemann zeta function. Expression (\ref{lo}) can be expanded into
\begin{equation}\label{b3bis}
\mathscr{I}=-\int_0^{1}\frac{\ln^2(t)\ln(1-t)}{t}\,dt-\int_0^{1}\frac{\ln^3(1-t)}{t}\,dt+2\int_0^{1}\frac{\ln (t)\ln^2(1-t)}{t}\,dt. 
\end{equation}
Lewin \cite{Lewin1981} provides the following expression (Eqs. (7.48) p. 199 and (7.61) p. 202):
\begin{equation}\label{Lewin761}
\mathrm{Li}_4(x)=\ln(x) \mathrm{Li}_3(x)-\frac{1}{2}\ln^2(x)\mathrm{Li}_2(x)-\frac{1}{2}\int_0^x\ln^2(t)\frac{\ln(1-t)}{t}\,dt.
\end{equation}
From Eq. (\ref{Lewin761}) one obtains, taking the limit $x\rightarrow 1$:
\begin{equation}\label{first}
\int_0^1\frac{\ln^2(t)\ln(1-t)}{1-t}\,dt=-2\zeta(4).
\end{equation}
Using Eq. (7.62) p. 203 of Ref. \cite{Lewin1981} or integrating further by parts in Eq (\ref{Lewin761}) leads to:
\begin{equation}\label{Lewin762}
\mathrm{Li}_4(x)=\ln(x) \mathrm{Li}_3(x)-\frac{1}{2}\ln^2(x)\mathrm{Li}_2(x)-\frac{1}{6}\ln^3(x)\ln(1-x)-\frac{1}{6}\int_0^x\frac{\ln^3(t)}{1-t}\,dt.
\end{equation}
Taking the limit $x\rightarrow 1$ in the latter equation and using the property that
\begin{equation}
\int_0^1\frac{\ln^3(t)}{1-t}\,dt=\int_0^1\frac{\ln^3(1-t)}{t}\,dt
\end{equation}
gives
\begin{equation}\label{second}
\int_0^x\frac{\ln^3(1-t)}{t}\,dt=-6\zeta(4).
\end{equation}
The last integral is the right-hand side of Eq. (\ref{b3bis}) is slightly more complicated. Let us start by differentiating $\ln^2(x)\ln^2(1-x)/2$:
\begin{equation}
\frac{1}{2}\frac{d}{dx}\ln^2(x)\ln^2(1-x)=\frac{\ln(x)}{x}\ln^2(1-x)-\ln^2(x)\frac{\ln(1-x)}{1-x}.
\end{equation}
Hence, integrating the latter expression implies
\begin{equation}\label{int1}
\frac{1}{2}\ln^2(x)\ln^2(1-x)=\int_0^x\frac{\ln(t)}{t}\ln^2(1-t)\,dt-\int_0^x\ln^2(t)\frac{\ln(1-t)}{1-t}\,dt.
\end{equation}
Now, let us change the variable from $x$ to $-x/(1-x)$ in Eq. (\ref{Lewin762}):
\begin{eqnarray}
\mathrm{Li_4}\left(\frac{-x}{1-x}\right)&=&\ln\left(\frac{x}{1-x}\right)\mathrm{Li_3}\left(\frac{-x}{1-x}\right)-\frac{1}{2}\ln^2\left(\frac{x}{1-x}\right)\mathrm{Li_2}\left(\frac{-x}{1-x}\right)\nonumber\\
& &+\frac{1}{6}\ln^3\left(\frac{x}{1-x}\right)\ln(1-x)\nonumber\\
& &+\frac{1}{6}\int_0^x\ln^3\left(\frac{t}{1-t}\right)\frac{dt}{1-t}.
\end{eqnarray}
Expanding the logarithm:
\begin{equation}
\int_0^x\left[\ln^3(t)-3\ln^2(t)\ln(1-t)+3\ln(t)\ln^2(1-t)-\ln^3(1-t)\right]\frac{dt}{1-t},
\end{equation}
we notice that, of the four terms of the latter expression, the first and third ones have already been dealt with. The fourth one is elementary and the second one can be evaluated in terms of the others. One finds, after some algebra (see Ref. \cite{Lewin1981}, formula (7.65), p. 204):
\begin{eqnarray}\label{int2}
\int_0^x\ln^2(t)\frac{\ln(1-t)}{1-t}\,dt&=&-2\left[\mathrm{Li}_4\left(-\frac{x}{1-x}\right)+\mathrm{Li}_4(x)-\mathrm{Li}_4(1-x)+\mathrm{Li}_4(1)\right]\nonumber\\
&=&+2\left[\ln(1-x)\mathrm{Li}_3(x)-\ln(x)\mathrm{Li}_3(1-x)\right]\nonumber\\
& &+2\ln(x)\ln(1-x)\mathrm{Li}_2(1-x)-\frac{\pi^2}{6}\ln^2(1-x)\nonumber\\
& &+\frac{1}{12}\ln^2(1-x)\left[6\ln^2(x)+4\ln(x)\ln(1-x)-\ln^2(1-x)\right]\nonumber\\
& &+2\mathrm{Li}_3(1)\ln\left(\frac{1}{1-x}\right).
\end{eqnarray}
As mentioned by Lewin, equating Eqs. (\ref{int1}) and (\ref{int2}) simplifies through cancellation of most of the terms due to a particular case of the inversion formula
\begin{equation}
    \mathrm{Li}_n(-x)+(-1)^n\mathrm{Li}_n\left(-\frac{1}{x}\right)=-\frac{1}{n!}\ln^n(x)+2\sum_{r=1}^{\lfloor\frac{n}{2}\rfloor}\frac{\ln^{n-2r}(x)}{(n-2r)!}\mathrm{Li}_{2r}(-1),
\end{equation}
where $\lfloor y\rfloor$ denotes the integer part of $y$, and one is then left with
\begin{equation}\label{third}
\int_0^{1}\frac{\ln t\ln^2(1-t)}{t}\,dt=-\frac{\zeta(4)}{2}.
\end{equation}
It is worth mentioning that an alternative proof of (\ref{third}) was provided by Connon \cite{Connon2008}. Inserting Eqs. (\ref{first}), (\ref{second}) and (\ref{third}) in the right-hand side of Eq. (\ref{b3bis}) gives the final result
\begin{equation}
\mathscr{I}=-[-2\zeta(4)]-[-6\zeta(4)]+2\left(-\frac{\zeta(4)}{2}\right)=7\zeta(4)=\frac{7\pi^4}{90}.
\end{equation}
One has therefore
\begin{empheq}[box=\fbox]{align}\label{final}
\zeta(4)=\frac{1}{7}\int_{-\infty}^{+\infty} z^2\, \frac{\ln(1+e^z)}{1+e^z}\,dz. 
\end{empheq}
It is easy to show that
\begin{equation}
\zeta(2)=\int_{-\infty}^{+\infty} \, \frac{\ln(1+e^z)}{1+e^z}\,dz. 
\end{equation}
The last two result may lead us to conjecture a general formula for $\zeta(2p+2)$ involving integrals of the kind
\begin{equation}
\int_{-\infty}^{+\infty} z^{2p}\, \frac{\ln(1+e^z)}{1+e^z}\,dz.
\end{equation}  
Using a computer algebra system \cite{Mathematica}, we obtained
\begin{equation}
\int_{-\infty}^{+\infty} z^{4}\, \frac{\ln(1+e^z)}{1+e^z}\,dz=\frac{279}{2} ~\zeta(6),
\end{equation}
as well as
\begin{equation}
\int_{-\infty}^{+\infty} z^{6}\, \frac{\ln(1+e^z)}{1+e^z}\,dz=5715~\zeta(8),
\end{equation}
or
\begin{equation}
\int_{-\infty}^{+\infty} z^{8}\, \frac{\ln(1+e^z)}{1+e^z}\,dz=\frac{804825}{2}~\zeta(10),
\end{equation}
but we could not find any general formula.

\section{Conclusion}

In this note, we proposed an integral representation for $\zeta(4)$. The corresponding expression is obtained using relations for polylogarithms. It is hoped that the present work will stimulate the derivation of further integral representations, and open the way to generalizations.

\section*{Acknowledgements}

I would like to thank Philippe Arnault for stimulating the present calculation in the course of a common work about statistical plasma physics.% \cite{Arnault2023}.


\begin{thebibliography}{99}

\bibitem{Tenenbaum2008} G. Tenenbaum, {\it Introduction à la th\'eorie analytique et probabiliste des nombres} (Belin, Paris, 2008) [in french].

\bibitem{Cartier1992} P. Cartier, {\it An introduction to zeta-functions}, in M. Waldschmidt, P. Moussa, J.-M. Luck et C. Itzykson (eds.), {\it From number theory to physics} (Springer-Verlag, Berlin, Heidelberg, 1992).

\bibitem{Borwein1995} D. Borwein and J. M. Borwein, On an intriguing integral and some series related to $\zeta(4)$, Proc. Amer. Math. Soc. {\bf 123}, 1191–1198 (1995).

\bibitem{Lewin1981} L. Lewin, {\it Polylogarithms and associated functions} (Elsevier North Holland, New York, Amsterdam, 1981).

\bibitem{Connon2008} D. F. Connon, {\it Some series and integrals involving the Riemann zeta function, binomial coefficients and the harmonic numbers. Volume I}, arXiv 0710.4022 (2008).\\
\url{https://arxiv.org/abs/0710.4022}

\bibitem{Mathematica} Wolfram Research, Inc., Mathematica, Version 13.2, Champaign, IL (2022).

%\bibitem{Arnault2023} P. Arnault, J. Racine, J.-P. Raucourt, A. Blanchet and J.-C. Pain, Sommerfeld expansion of electronic entropy in an inferno-like average atom model, Phys. Rev. B {\bf 108}, 085115 (2023).

\end{thebibliography}
\end{document}